\documentclass{amsart}
\usepackage{comment}

\usepackage[all]{xy}
\usepackage{pgf}
\usepackage{tikz}
\usetikzlibrary{arrows,automata,snakes}
\usepackage[utf8]{inputenc}
\usepackage{amsmath,amsfonts,amsthm,amssymb}
\usepackage{xfrac}
\usepackage{verbatim}
\usepackage[english]{babel}
\usepackage{todonotes}
\usepackage[backref]{hyperref}
\usepackage{cite}

\DeclareMathOperator{\Aut}{Aut}

\DeclareMathOperator{\GL}{GL}

\DeclareMathOperator{\divisor}{div}

\newcommand{\Q}{{\mathbb Q}}
\newcommand{\Z}{{\mathbb Z}}
\newcommand{\C}{{\mathbb C}}

\newcommand{\F}{{\mathbb F}}
\newcommand{\PP}{{\mathbb P}}

\newcommand{\cL}{\mathcal{L}}

\newcommand{\cN}{\mathcal{N}}

\newcommand{\ff}{\mathfrak{f}}

\newcommand{\set}[1]{\left\lbrace #1 \right\rbrace}
\newcommand{\field}[1]{\mathbb{#1}}
\renewcommand{\P}{\field{P}}

\begin {document}

\newtheorem{thm}{Theorem}
\newtheorem{lem}{Lemma}[section]
\newtheorem{prop}[lem]{Proposition}

\theoremstyle{definition}

\theoremstyle{remark}

\title[]{Elliptic Curves over Totally Real Cubic Fields\\
are Modular}

\author{Maarten Derickx}
\address{Department of Mathematics, MIT 2-252B, 77 Massachusetts Ave, Cambridge, MA 02139}
\email{drx@mit.edu}

\author{Filip Najman}

\address{Department of Mathematics, Faculty of Science, University of Zagreb\\
Bijeni\v cka cesta 30\\
10000 Zagreb\\
Croatia}
\email{fnajman@math.hr }

\author{Samir Siksek}

\address{Mathematics Institute\\
	University of Warwick\\
	CV4 7AL \\
	United Kingdom}

\email{s.siksek@warwick.ac.uk}

\date{\today}
\thanks{
Derickx is supported by Simons Foundation grant 550033. 
Najman is supported by the QuantiXLie Centre of Excellence, a project
co-financed by the Croatian Government and European Union through the
European Regional Development Fund - the Competitiveness and Cohesion
Operational Programme (Grant KK.01.1.1.01.0004). Siksek is supported by
EPSRC {\em LMF: L-Functions and Modular Forms} Programme Grant
EP/K034383/1.
}

\begin{abstract}
We prove that all elliptic curves defined over totally real
cubic fields are modular. This builds on previous
work of Freitas, Le Hung and Siksek, who proved modularity of
elliptic curves over real quadratic fields, as well as  
recent breakthroughs due to Thorne and to Kalyanswamy.
\end{abstract}
\maketitle
%\makeatletter
%\renewcommand{\tocsection}[3]{%
%  \indentlabel{\@ifnotempty{#2}{\ignorespaces#1 #2.\quad}}#3\dotfill
%}
%\let\tocsubsection\tocsection
%\makeatother
%\tableofcontents

\section{Introduction}

Let $K$ be a totally real number field and let $E$ be an elliptic curve
over $K$ with conductor $\cN$. It is conjectured that such a curve $E$ is \textbf{modular} in the following
sense: there is a level $\cN$ Hilbert newform $\ff$ over $K$ of parallel weight $2$
and rational Hecke eigenvalues such that
$L(E,s)=L(\ff,s)$, where the $L$-function on the left is the Hasse--Weil $L$-function of $E$,
and the $L$-function on the right is the Hecke $L$-function of $\ff$. This \textbf{modularity
conjecture} is the natural generalization to totally real fields
of the Shimura--Taniyama conjecture for elliptic curves over the rationals.
The latter is a celebrated theorem due to Wiles \cite{Wiles}, Breuil, Conrad, Diamond
and Taylor \cite{modularity}. The earliest results towards the modularity conjecture
for elliptic curves
going beyond the rationals were due to Jarvis and Manoharmayum \cite{JM},
and established modularity of semistable elliptic
curves over $\Q(\sqrt{2})$ and $\Q(\sqrt{17})$. In the last 10 years there has
been a dramatic strengthening of modularity lifting theorems
due to, for example,
Breuil and Diamond \cite{BreuilDiamond},
Kisin \cite{Kisin}, Gee \cite{Gee}, and
Barnet-Lamb, Gee and Geraghty \cite{BGG1}, \cite{BGG2}.
The following theorem is easily deduced from these aforementioned modularity lifting theorems
and by now standard
modularity switching arguments due to Wiles and to Manoharmayum \cite{Mano1}.
(A proof is given by Freitas, Le Hung and Siksek \cite[Theorems 3 and 4]{FHS}
but the arguments are well-known).
\begin{thm}\label{thm:357}
Let $p=3$, $5$ or $7$. Let $E$ be an elliptic curve over a totally real field $K$, and write $\overline{\rho}_{E,p}$ for its mod $p$ representation.
Suppose that $\overline{\rho}_{E,p}(G_{K(\zeta_p)})$ is absolutely irreducible.
Then $E$ is modular.
\end{thm}
A hypothetical non-modular $E/K$ would therefore necessarily have small mod $p$ image
for $p=3$, $5$, $7$ and would give rise to a $K$-point on one of a number of modular
curves---we make this precise later. In \cite{FHS}, the real quadratic points
of these modular curves are shown to be either cuspidal,
or to correspond to elliptic curves that have complex multiplication,
or rational $j$-invariants,
or that are $\Q$-curves. The authors deduce the following.
\begin{thm}[Freitas, Le Hung and Siksek]\label{thm:FHS}
Elliptic curves over real quadratic fields are modular.
\end{thm}
Whilst direct computation is used in \cite{FHS} to  study the real quadratic
points on the relevant modular curves, Thorne \cite{Thorne5} uses Iwasawa
theory to control points on some of these modular
curves over $\Z_p$-extensions of $\Q$,
and deduces the following remarkable theorem.
\begin{thm}[Thorne]
Let $p$ be a prime, and let $K$ be a number field which is contained in the cyclotomic $\Z_p$-extension
of $\Q$. Let E be an elliptic curve over $K$. Then $E$ is modular.
\end{thm}

Recently Theorem~\ref{thm:357} was
substantially strengthened in the cases $p=5$ and $p=7$, respectively
by Thorne \cite{Thorne} and Kalyanswamy \cite{Kalyanswamy}.
This means that several difficult steps in the proof of Theorem~\ref{thm:FHS}
can now be eliminated.
In this paper we build on these theorems of Thorne and Kalyanswamy
to prove the following.
\begin{thm}\label{thm:main}
Let $K$ be a totally real cubic number field. Let $E$ be an elliptic curve over $K$.
Then $E$ is modular.
\end{thm}

\bigskip

The computations in this paper were carried out
in the computer algebra system \texttt{Magma} \cite{MAGMA}.
The reader can find the \texttt{Magma} scripts
for verifying these 
computations
at: \newline
\url{http://homepages.warwick.ac.uk/staff/S.Siksek/progs/cubicmodularity/}

\section{Images mod $3$, $5$, $7$ and modularity}
Let $p\ge 3$ be a prime.
Write $B(p)$ for a Borel subgroup of $\GL_2(\F_p)$,
and $C_\mathrm{s}(p)$ and $C_\mathrm{ns}(p)$ respectively
for a split and non-split Cartan subgroups.
Let $C_\mathrm{s}^+(p)$ and $C_\mathrm{ns}^+(p)$
respectively be their normalizers.
The three groups
$B(p)$, $C_\mathrm{s}^+(p)$ and $C_\mathrm{ns}^+(p)$
correspond to modular curves defined over $\Q$
which are usually denoted by
$X_0(p)$, $X_\mathrm{split}(p)$ and $X_\mathrm{nonsplit}(p)$.
Instead we shall mostly follow the notation of \cite{FHS}
and denote these modular curves
by: $X(\mathrm{b}p)$, $X(\mathrm{s}p)$ and $X(\mathrm{ns}p)$.
The following is well-known. For a proof, see \cite[Section 2.3]{FHS}.
\begin{prop}\label{prop:large}
Let $E$ be an elliptic curve over a totally real field $K$,
and let $p \ge 3$ be
a rational prime.
%and write $\overline{\rho}=\overline{\rho}_{E,p}$.
%Then
%\begin{enumerate}
%\item[(i)] $\overline{\rho}(G_{K(\zeta_p)})=\overline{\rho}(G_K) \cap \SL_2(\F_p)$.
%\item[(ii)]
If $\overline{\rho}(G_{K(\zeta_p)})$ is absolutely reducible,
then $\overline{\rho}(G_K)$ is contained either in a Borel subgroup,
or in the normalizer of a Cartan subgroup. In this case $E$ gives
rise to a non-cuspidal $K$-point on $X(\mathrm{b}p)$, $X(\mathrm{s}p)$
or $X(\mathrm{ns}p)$.
%\end{enumerate}
\end{prop}

When the additional assumption $K \cap \Q(\zeta_p)=\Q$ is satisfied, we can substantially improve on
Proposition~\ref{prop:large}. We shall only need
this improvement for $p=3$; in this case, as $K$ is real, it is certainly
satisfies the additional assumption. The following is
part of Proposition 4.1
in \cite{FHS}.
\begin{prop}\label{prop:large2}
Let $K$ be a totally real number field
and let $E$ be an elliptic curve over $K$.
%Let $p=3$, $7$ be a rational prime and Suppose
%\[
%\]
Write
$\overline{\rho}=\overline{\rho}_{E,3}$.
Suppose %$\overline{\rho}$ is irreducible but
$\overline{\rho}(G_{K(\zeta_p)})$ is absolutely reducible.
Then $\overline{\rho}(G_K)$ is conjugate to
a subgroup of $B(3)$ or $C_\mathrm{s}^+(3)$.
\end{prop}

We shall need the following strengthening of Theorem~\ref{thm:357}
for $p=5$ due to Thorne \cite{Thorne}.
\begin{thm}[Thorne] \label{thm:Thorne}
Let $K$ be a totally real field and $E$ an elliptic curve over $K$.
Suppose $5$ is not a square in $K$, and $\overline{\rho}_{E,5}$
is irreducible. Then $E$ is modular.
\end{thm}

We shall need the following theorem of
Kalyanswamy \cite[Proposition 4.3 and Theorem 4.4]{Kalyanswamy}
which improves on Theorem~\ref{thm:357} for $p=7$.
\begin{thm}[Kalyanswamy]\label{thm:Kalyanswamy}
Let $K$ be a totally real field and $E$ an elliptic curve over $K$.
Suppose
\begin{itemize}
\item $K \cap \Q(\zeta_7)=\Q$.
\item $\overline{\rho}_{E,7}$ is irreducible.
\item $\overline{\rho}_{E,7}(G_K)$ is not conjugate to
a subgroup of $C_\mathrm{ns}^+(7)$.
%the group
%$H_2$.
%\[
%H_2=
%\Biggl\langle
%\begin{pmatrix}
%0 & 5\\
%3 & 0
%\end{pmatrix},
%\begin{pmatrix}
%5 & 0\\
%3 & 2
%\end{pmatrix}
%\Biggr\rangle \subset \GL_2(\F_7).
%\]
\end{itemize}
Then $E$ is modular.
\end{thm}
Kalyanswamy's theorem is somewhat more precise, but we shall not
need its full strength.

In this paper we shall deal with four modular curves, that
we denote by $X(\mathrm{b}3,\mathrm{b}5)$, $X(\mathrm{s}3,\mathrm{b}5)$,
$X(\mathrm{b}5,\mathrm{b}7)$, $X(\mathrm{b}5,\mathrm{ns}7)$.
For the details of the notation, and the precise modular interpretation
of these curves we refer to \cite[Section 2.2.2]{FHS}.
These are often conveniently thought of as normalizations
of fibre products: if $p$, $q$
are distinct primes and $\mathrm{u}$, $\mathrm{v} \in \{\mathrm{b},\mathrm{s},
\mathrm{ns}\}$ then $X(\mathrm{u}p,\mathrm{v}q)$
is the normalization of $X(\mathrm{u}p) \times_{X(1)} X(\mathrm{v}q)$.

\section{Modularity of Elliptic Curves over $\Q(\zeta_7)^+$}
In this section prove Theorem~\ref{thm:main} for $K=\Q(\zeta_7)^+$.
\begin{lem}\label{lem:qz7}
Let $K=\Q(\zeta_7)^+$. Let $E$ be an elliptic curve defined
over $K$. Then $E$ is modular.
\end{lem}
\begin{proof}
By Theorem~\ref{thm:Thorne} we may suppose that $\overline{\rho}_{E,5}$
is reducible. By Theorem~\ref{thm:357} and Proposition~\ref{prop:large2}
we may suppose that the image of $\overline{\rho}_{E,3}$
is contained in $B(3)$ or $C_{\mathrm{s}}^+(3)$.
Thus $E$ gives rise to a non-cuspidal $K$-point on
one of the two modular curves $X(\mathrm{b}3,\mathrm{b}5)$,
$X(\mathrm{s}3,\mathrm{b}5)$.
%\[
%X(\mathrm{b}3,\mathrm{b}5):=X(\mathrm{b}3) \times_{X(1)} X(\mathrm{b}5),
%\qquad
%X(\mathrm{s}3,\mathrm{b}5):=X(\mathrm{s}3) \times_{X(1)} X(\mathrm{b}5).
%\]
It is shown in \cite[Section 5.4.2]{FHS}
that these are in fact elliptic curves defined over $\Q$
with Cremona labels \texttt{15A1} and \texttt{15A3}.
We computed the Mordell--Weil groups $X(K)$ for
$X=X(\mathrm{b}3,\mathrm{b}5)$,
$X(\mathrm{s}3,\mathrm{b}5)$ using \texttt{Magma}. In both cases
we found
\[
X(K)=X(\Q) \cong \Z/4\Z \oplus \Z/2\Z.
\]
In particular $E$ gives rise to $\Q$-point on $X$
and so is a twist of an elliptic curve defined over $\Q$.
It is therefore modular by \cite{modularity}.
\end{proof}

\section{An overview of the proof of Theorem~\ref{thm:main}}
The remainder of the paper is devoted to the proof of
the following two theorems.
\begin{thm}\label{thm:b5b7}
Let $K$ be a totally real cubic field.
Then $X(\mathrm{b}5,\mathrm{b}7)(K)$
consists only of cusps.
%\[
%X=X(\mathrm{b}5,\mathrm{b}7)
%:=X(\mathrm{b}5) \times_{X(1)} X(\mathrm{b}7).
%\]
%Then $X(K)$ consists only of cusps.
\end{thm}

\begin{thm}\label{thm:cubic_pts}
Let $K$ be a cubic field.
Then $X(\mathrm{b}5,\mathrm{ns}7)(K)$ consists
only of cusps.
%\begin{equation}\label{eqn:fibprod}
%X=X(\mathrm{b}5,\mathrm{ns}7):=X(\mathrm{b}5) \times_{X(1)} X(\mathrm{ns}7).
%\end{equation}
%Let $K$ be a cubic field. Then $X(K)$ consists only of cusps.
\end{thm}
\begin{comment}
\begin{thm}\label{thm:cubic_pts}
Let
\[
X(\mathrm{b}3,\mathrm{b}5,\mathrm{ns}7):=X(\mathrm{b}3,\mathrm{b}5) \times_{X(1)} X(\mathrm{ns}7),
\]
and
\[
X(\mathrm{s}3,\mathrm{b}5,\mathrm{ns}7):=X(\mathrm{s}3,\mathrm{b}5) \times_{X(1)} X(\mathrm{ns}7).
\]
The curves $X(\mathrm{b}3,\mathrm{b}5,\mathrm{ns}7)$ and $X(\mathrm{s}3,\mathrm{b}5,\mathrm{ns}7)$ have no noncuspidal cubic points.
\end{thm}
\end{comment}
In this section we complete the proof of Theorem~\ref{thm:main}
assuming Theorems~\ref{thm:b5b7} and~\ref{thm:cubic_pts}.
Let $K$ be a totally real cubic field and $E$ an elliptic
curve over $K$. We would like to show that $E$ is modular.
Suppose that it is not modular.
%We start by summarising what we can assume about the
%mod $5$ and mod $7$ images of $E$.
%\begin{enumerate}
%\item[(a)] By Theorem~\ref{thm:357} and Proposition~\ref{prop:large2}
%we may assume that the image of $\overline{\rho}_{E,3}$
%is contained in either $B(3)$ or $C_{\mathrm{s}}^+(3)$
%and thus $E$ gives rise to a non-cuspidal $K$-point
%on $X(\mathrm{b}3)$ or $X(\mathrm{s}3)$.
%\item[(a)]
By Theorem~\ref{thm:Thorne},
the representation $\overline{\rho}_{E,5}$ is reducible.
%
%and so $E$ gives rise to a non-cuspidal $K$-point
%on $X(\mathrm{b}5)$.
%\item[(b)]
By Lemma~\ref{lem:qz7} we know $K \ne \Q(\zeta_7)^+$
and thus $K \cap \Q(\zeta_7)=\Q$.
%By Lemma~\ref{lem:b5b7}
%we may suppose that $\overline{\rho}_{E,7}$ is irreducible.
Thus by Theorem~\ref{thm:Kalyanswamy}
the image of $\overline{\rho}_{E,7}$ is contained in $B(7)$
or in
$C_{\mathrm{ns}}^+(7)$.
% and so $E$ gives rise to a non-cuspidal $K$-point
%on either $X(\mathrm{b}7)$ or on $X(\mathrm{ns}7)$.
%\end{enumerate}
Thus if $E$ is not modular then $E$ gives
rise to a non-cuspidal $K$-point on either $X(\mathrm{b}5,\mathrm{b}7)$
or $X(\mathrm{b}5,\mathrm{ns}7)$.
%Suppose first that the image of $\overline{\rho}_{E,7}$ is contained in $B(7)$.
%From (a) we see that $E$ gives rise to a non-cuspidal $K$-point on
%$X(\mathrm{b}5,\mathrm{b}7)$.
But by Theorems~\ref{thm:b5b7} and~\ref{thm:cubic_pts} there are
no such points, giving a contradiction.
%This contradicts
%Theorem~\ref{thm:b5b7}. We conclude that $\overline{\rho}_{E,7}$
%is contained in $C_{\mathrm{ns}}^+(7)$. Now from (a), (b) we see that $E$
%gives rise to a non-cuspidal
%$K$-point on $X(\mathrm{b}5,\mathrm{ns}7)$.
%This contradicts
%Theorem~\ref{thm:cubic_pts}.

In summary, to prove Theorem~\ref{thm:main} all we
have to do is prove Theorems~\ref{thm:b5b7} and~\ref{thm:cubic_pts}.

\section{Proof of Theorem~\ref{thm:b5b7}}\label{sec:b5b7}
Let $X=X(\mathrm{b}5,\mathrm{b}7)$ (in standard notation
denoted by $X_0(35)$).
It is known that $X$ has four $\Q$-points and that these
are cusps.
Let $K$ be a totally real cubic field.
For the proof of Theorem~\ref{thm:b5b7}
it will be sufficient to show that $X(K)=X(\Q)$.
Suppose $P \in X(K) \backslash X(\Q)$.
Let $P_1$, $P_2$, $P_3$
be the conjugates of $P$ given by the three embeddings of $K$
in $\overline{\Q}$, and write $D=P_1+P_2+P_3$.
Then $D$ is an irreducible $\Q$-rational divisor on $X$
of degree $3$. We shall determine all the irreducible
$\Q$-rational divisors of degree $3$ on $X$
and show that none of them arise from totally real cubic
points, giving a contradiction.

The arithmetic of $X$ and its Jacobian are studied
in \cite[Section 5.1]{FHS}.
The curve $X$ is hyperelliptic of genus $3$. A model for $X$,
derived by Galbraith \cite[Section 4.4]{Galbraith}, is given by
%Galbraith \cite[Section 4.4]{Galbraith} derives the following model for $X$,
\begin{equation}\label{eqn:35}
X \; : \; y^2=(x^2 + x - 1)(x^6 - 5 x^5 - 9 x^3 - 5 x - 1).
\end{equation}
Write $\infty_{\pm}$ for the two points at infinity.
%There are precisely four $\Q$-points on $X$ and these
%are the four cusps: $\infty_{\pm}$, $(0,\pm 1)$.
Write $J$ for $J_0(35)$---the Jacobian of $X$.
% Using $2$-descent
%\cite{Stoll} we find that $J(\Q)$ has Mordell-Weil rank $0$.
%As $J$ has good reduction at $3$, the map $J(\Q) \rightarrow J(\F_3)$
%is injective. The group generated by differences of the four obvious
%rational points $(0,\pm 1)$, $\infty_{\pm}$ (where the latter
%are the points at infinity) surjects onto
%$J(\F_3) \cong \Z/24\Z \times \Z/2\Z$, and is therefore equal
%to $J(\Q)$.
%We find
Then
\[
J(\Q)= \frac{\Z}{24 \Z} \cdot \left[ \infty_{-}-\infty_{+} \right]
+\frac{\Z}{2 \Z} \cdot \left[3(0,-1) -3 \infty_{+} \right].
\]
Let $D_1,\dots,D_{48}$ be rational divisors of degree $0$
on $X$ representing the $48$ classes in $J(\Q)$,
and let $D_i^\prime=D_i+3 \infty_{+}$. Recall
that $D$ is an irreducible $\Q$-rational divisor of degree $3$.
Then $D \sim D_i^\prime$ for some $i$.
We shall write $\cL(D_i^\prime)$ for the Riemann--Roch space
corresponding to $D_i^\prime$ and $\lvert D_i^\prime \rvert$ for
the corresponding complete linear system. By Riemann--Roch and
Clifford's inequality, $\dim \cL(D_i^\prime) =1$ or $2$.
Moreover, if $\dim \cL(D_i^\prime)=2$, then $\lvert D_i^\prime \rvert$
contains a base point
(c.f.\ \cite[Chapter I, Exercise D.9]{ACGH}), and therefore cannot
contain an irreducible divisor. To sum up, $D \sim D_i^\prime$
for some $1 \le i \le 48$ such that $\dim \cL(D_i^\prime)=1$.
We computed these spaces using \texttt{Magma}; for this \texttt{Magma}
uses an algorithm of Hess \cite{Hess}. We found
that $\dim \cL(D_i^\prime)=1$ for precisely $44$ of the
$48$ divisors $D_i^\prime$. For these, letting $f_i$ be
a $\Q$-basis for $\cL(D_i^\prime)$, gives $D=D_i^\prime+\divisor(f_i)$
for some $i$. We found that precisely $28$ of the
effective degree $3$ divisors $D_i^\prime+\divisor(f_i)$ are
irreducible. However, all of these split over a cubic field
with a complex embedding giving the required contradiction.

\section{The modular curve $X(\mathrm{b}5,\mathrm{ns}7)$}
\begin{comment}
In the Section~\ref{sec:rat_pts} we show that the
cubic points on the curves
$X(\mathrm{b}3,\mathrm{b}5,\mathrm{ns}7)$,
$X(\mathrm{s}3,\mathrm{b}5,\mathrm{ns}7)$
are cuspidal.
%From \cite{FHS}, these two curves have genera $153$ and $73$
%respectively.
It is somewhat challenging to work with these
curves directly and in this section we first study
a common quotient of both curves:
\begin{equation}\label{eqn:fibprod}
X(\mathrm{b}5,\mathrm{ns}7):=X(\mathrm{b}5) \times_{X(1)} X(\mathrm{ns}7).
\end{equation}
\end{comment}
We shall henceforth restrict our attention to
$X(\mathrm{b}5,\mathrm{ns}7)$. To simplify the notation
we write $X=X(\mathrm{b}5,\mathrm{ns}7)$. We denote
the Jacobian of $X$ by $J=J(\mathrm{b}5,\mathrm{ns}7)$.
%This is a common quotient of the aforementioned two curves
%as $H_2$ is contained in $C_{\mathrm{ns}}^+(7)$.
The curve $X$
and its Jacobian $J$ are studied
in Le Hung's thesis \cite[Section 6.4]{lehung} and we make
extensive use of his results. In particular, this curve is non-hyperelliptic
and has
genus $6$.
\subsection{The Jacobian $J=J(\mathrm{b}5,\mathrm{ns}7)$}
Le Hung shows that
\[
J \sim A_1 \times A_2 \times A_3
\]
where $\sim$ here denotes isogeny over $\Q$, and $A_1$, $A_2$, $A_3$
are modular abelian surfaces defined over $\Q$. Morever the $A_i$ are absolutely simple.
The involution $w_5$ on $J$
is compatible with the isogeny and
acts by multiplication by $1$, $-1$, $-1$
respectively on $A_1$, $A_2$, $A_3$. The analytic ranks of $A_1$, $A_2$, $A_3$
are respectively $2$, $0$, $0$. In particular, by the work
of Kolyvagin and Logach\"{e}v \cite{KL}, the Mordell--Weil groups
$A_2(\Q)$ and $A_3(\Q)$ are torsion.
We immediately deduce the following.
\begin{lem}\label{lem:Ators}
Let $A/\Q$ be the abelian subvariety of $J$ that is the image
of $w_5-1$.
Then $A \sim A_2 \times A_3$ has dimension $4$. Moreover,
the Mordell--Weil group $A(\Q)$ is torsion.
\end{lem}
%either
%trivial, or isomorphic to $\Z/7\Z$. \margnote{The way this lemma
%was originally worded would require us to prove that there
%is really a $\Z/7\Z$ subgroup.}
%$J(\mathrm{b}5,\mathrm{ns}7)(\Q) \cong \Z^2\times \Z/7\Z$,
%and furthermore this isomorphism can
%be chosen so that $w_5$ acts as the identity on the $\Z^2$ factor and as $-1$
%on the $\Z/7\Z$ factor.
%In particular if $A \subseteq J(\mathrm{b}5,\mathrm{ns}7)$ is the abelian
%subvariety that is the image of
%$w_5-1 : J(\mathrm{b}5,\mathrm{ns}7) \to J(b\mathrm{5},\mathrm{ns}7)$
%then $A(\Q) =
%A(\Q)_{tors} \cong \Z/7\Z$.
%\begin{proof}
%The factor $A_1$ is the only one of the $A_i$ with (possibly) positive
%rank. As $w_5$ is the identity on $A_1$, the endomorphism
%$w_5-1$ kills $A_1$. Thus
%$(w_5-1) J(\mathrm{b}5,\mathrm{ns}7)(\Q)$ is torsion.
%\end{proof}

\subsection{Le Hung's model for $X=X(\mathrm{b}5,\mathrm{ns}7)$}
We need a good model for $X(\mathrm{b}5,\mathrm{ns}7)$.
Le Hung \cite[p. 47]{lehung} gives a
model which will be a good starting point for us.
We briefly sketch Le Hung's derivation of his model,
but work with projective rather than affine coordinates.
Later we explain how to derive a better model.
The curves $X(\mathrm{b}5)$ and $X(\mathrm{ns}7)$
are both isomorphic to $\PP^1$ over $\Q$.
Let
\[
F_1(x_1,x_2)=(x_1^2+10x_1 x_2+5x_2^2)^3 , \qquad F_2:=x_1 x_2^5,
\]
\[
G_1(y_1,y_2)=64 \cdot \left(y_1\cdot (y_1^2+7y_2^2) \cdot (y_1^2-7y_1 y_2+14y_2^2)\cdot (5y_1^2-14y_1 y_2-7y_2^2)\right)^3 ,
\]
and
\[
G_2(y_1,y_2)=(y_1^3-7y_1^2 y_2+7y_1 y_2^2+7y_2^3)^7.
\]
For appropriate choices of projective coordinates $(x_1:x_2)$
for $X(\mathrm{b}5)$ and $(y_1:y_2)$ on $X(\mathrm{ns}7)$,
the $j$-maps are given by %(\cite[p. 47]{lehung})
\[
j : X(\mathrm{b}5) \rightarrow X(1), \qquad (x_1 : x_2) \mapsto (F_1(x_1,x_2) : F_2(x_1,x_2)),
\]
and
\[
j: X(\mathrm{ns}7) \rightarrow X(1) ,
\qquad (y_1,y_2) \mapsto (G_1(y_1,y_2) : G_2(y_1,y_2)).
\]
As $X$ is the normalization of $X(\mathrm{b}5)\times_{X(1)} X(\mathrm{ns}7)$
 we immediately deduce the following
model for $X$ in $\PP^1 \times \PP^1$:
\[
C \; : \qquad F_1(x_1,x_2) G_2(y_1,y_2)=F_2(x_1,x_2)G_1(y_1,y_2).
\]
The curve $X$ is the normalization of this model.
%
% and in Section~\ref{sec:} we use this and the previous results
%to complete the proof of Theorem~\ref{thm:main}.
%One computational issue that we ran into when trying to do the necessary
%computations is that we needed explicit equations for the involution $w_5$ on
%$X(\mathrm{b}5,\mathrm{ns}7)$ for the computations in Section \ref{sec:rat_pts}. In \cite[p.
%47]{lehung} there is already the nice equation
%\begin{equation}\label{eqn:lehung}
%\[
%C \quad : \quad \frac{(x^2+10x+5)^3}{x} =
%\frac{ 64(\phi(\phi^2+7)(\phi^2-7\phi+14)(5\phi^2-14\phi-7))^3}{(\phi^3-7\phi^2+7\phi+7)^7};
%\]
%\end{equation}
%To understand the cusps of $X(\mathrm{b}5,\mathrm{ns}7)$
%it is more convenient to work with Le Hung's model \eqref{eqn:lehung}.
 The parameterization $(x_1 : x_2)$ on $X(\mathrm{b}5)$ is chosen
so that the $0$ and $\infty$ cusps are $(x_1:x_2)=(0:1)$ and $(x_1:x_2)=(1:0)$
respectively. We shall denote these by $a_0$, $a_\infty$.
Let $\zeta_7$ be a primitive $7$-th root of unity.
Let $\eta=2(\zeta_7^3+\zeta_7^{-3})+3 \in \Q(\zeta_7)^+$.
Then $F_2(\eta : 1)=0$. The three cusps of $X(\mathrm{ns}7)$
are $(\eta:1)$ and its Galois conjugates.
It follows that the cusps of $X$
are the points belonging to the normalization of $C$ lying above the points
$(x_1:x_2,y_1:y_2)=(0:1, \eta:1)$, $(1:0,\eta:1)$
and their Galois conjugates. Although these points on $C$ are singular,
it is easy to check (c.f. \cite[Section 5.5.1]{FHS})
that there is only one point on the normalization above each, and to deduce:
\begin{itemize}
\item $X$
has two Galois orbits of cusps, both of degree $3$ and defined
over $\Q(\zeta_7)^+$, which we denote by $c_0$, $c_\infty$;
\item The three cusps in $c_0$ map to
$a_0$, and the three cusps in $c_\infty$ map to $a_\infty$ on $X(\mathrm{b}5)$.
\item The divisor of $x_1/x_2$ interpreted as a function on
$X$ is
$7\cdot \left( c_0-c_\infty \right)$.
 In particular,
the class
$\left[c_0-c_\infty \right]$
 is an element of order $1$ or $7$. There are several ways to show
that the divisor $c_0-c_\infty$ is not principal, and
so its class has order $7$. One way is by direct computation
using \texttt{Magma}, working with the model $D$ introduced below.
Here is another way: we shall show below that $X$ has gonality $4$.
As $c_0$, $c_\infty$ have degree $3$ they cannot be linearly equivalent. 
\end{itemize}
%for $X(\mathrm{b}5,\mathrm{ns}7)$ and on this model $w_5^*(x)=\frac {125} x$. However, a formula
%for $w_5^*(\phi)$ is not given and asking \texttt{Magma} to compute the
%automorphism group of the curve $X(b5,ns7)$ is a computation that did not
%terminate after several days.

\subsection{A plane degree $6$ model for $X=X(\mathrm{b}5,\mathrm{ns}7)$}
We used \texttt{Magma} to compute,
starting with the model $C$,
the canonical map and its image. The latter is indeed
a smooth genus $6$ curve cut out in $\PP^5$ by six homogeneous degree $2$
polynomials. By the Enriques--Babbage Theorem \cite[p. 124]{ACGH},
we know that $X$ is neither trigonal,
nor isomorphic to a plane quintic. Moreover, as the factors
$A_i$ of the Jacobian are $2$-dimensional and absolutely simple,
we see that the curve is not bi-elliptic.
%and so in particular it has gonality $4$.
It follows (c.f. \cite[209--210]{ACGH}) that $X$
has gonality $4$ and
a degree $6$ planar model, with four ordinary double points as
singularities. We used the inbuilt \texttt{Magma} implementation
for writing down this model,
%However the \texttt{Magma} command $Genus6PlaneCurveModel$ happily computes a degree 6 planar model of the curve. It turns out that this degree $6$ model has 4 ordinary double point singularities,
and found that two of the four double points are defined over $\Q(i)$
and the other two over $\Q(\sqrt{5})$.
After applying a $\Q$-rational automorphism of $\P^2$ to slightly simplify
this degree $6$ model, it is given by the following equation:
\begin{align}
D\; : \; 5u^6 - 50u^5v + 206u^4v^2 - 408u^3v^3 + 321u^2v^4 + 10uv^5 - 100v^6 + 9u^4w^2 - \nonumber\\
%\label{eqn:uvw} \\
60u^3vw^2 + 80u^2v^2w^2 + 48uv^3w^2 + 15v^4w^2 + 3u^2w^4 - 10uvw^4 + 6v^2w^4 - w^6=0. \nonumber
\end{align}
On this model $D$ the double points are
%$(\pm i:0:1)$ and
%$(0:\pm \frac 1 {\sqrt 5}:1)$.
\[
p_1 = (i:0:1),\quad p_2 = (-i:0:1), \quad p_3= (0: \frac 1 {\sqrt 5}:1), \quad p_4 =
(0:-\frac 1 {\sqrt 5}:1).
\]

It is clear that $D$ has an automorphism $(u:v:w)\mapsto (-u:-v:w)$.
The curve $X$ has an obvious modular
involution which is $w_5$. The following lemma
proves that $w_5$ coincides with the automorphism $(u:v:w)\mapsto (-u:-v:w)$.
%This automorphism actually has to be $w_5$ because of the following lemma.
\begin{lem}
	The $\Q$-rational automorphism group of $X(b5,ns7)$ is generated by $w_5$, i.e. $\Aut_\Q (X) = \langle w_5 \rangle \cong \Z/2\Z$.
\end{lem}
\begin{proof}
%The singular points on the model $D$ are
% $p_1 = (i:0:1), p_2 = (-i:0:1), p_3= (0: \frac 1 {\sqrt 5}:1)$ and $p_4 =
%(0:-\frac 1 {\sqrt 5}:1)$.
As described in \cite[p 210--211.]{ACGH} a degree $6$ planar curve with
four
ordinary double points such as $D$ %(which $X(b5,ns7)$ is)
has exactly five different
$g_4^1$. Namely, one given by the pencil of quadrics going through all four
points, and the other four coming from the pencil of lines through each of the
$p_i$. Since none of the $p_i$ are $\Q$-rational, only the first $g_4^1$ is
defined over $\Q$. Now every $g_6^2$ on
%$X(b5,ns7)$
such a curve is residual to a $g_4^1$. This means that there is only one $\Q$-rational $g_6^2$, namely the one
corresponding to the degree $6$ model given by $u$, $w$, $v$ above.
In particular
every $\Q$-rational automorphism has to come from an automorphism $h: \P^2_\Q
\to \P^2_\Q$ in the degree $6$ model.
Such an automorphism $h$ has to preserve
the singular locus $\set{p_1,p_2,p_3,p_4}$ and is in fact uniquely determined
by what it does on this singular locus. Of the $24$ automorphisms of
$\P^2_{\overline{\Q}}$ preserving $\set{p_1,p_2,p_3,p_4}$, only the ones of the
form $(u:v:w) \mapsto (\pm u: \pm v: w)$ are $\Q$-rational. One easily sees
that of these four only the identity and  $(u:v:w)\mapsto (-u:-v:w)$ are
actually automorphisms of the curve.
\end{proof}

Transferring $c_0$ and $c_\infty$ to our new model $D$, we find that they
respectively are
the Galois orbits of the following
two points defined over $\Q(\eta)=\Q(\zeta_7)^+$ by
\[
(-4\eta^2 + 21\eta + 7: -\eta^2 + 7\eta : 14), \qquad
(4\eta^2 - 21\eta - 7: \eta^2 - 7\eta : 14).
\]
We note that these are interchanged by $w_5 : (u:v:w)\mapsto (-u:-v:w)$
as expected.

\subsection{The Mordell--Weil group $A(\Q)$}
In Lemma~\ref{lem:Ators} we defined the abelian subvariety
$A$ of $J$ as the image of $w_5-1$
and observed that $A(\Q)$ is torsion. We can now pin down
$A(\Q)$ precisely. In particular, applying the function field
class group algorithm of Hess \cite{Hess} (implemented in
\texttt{Magma}) to our model $D$ obtain
\[
J(\F_3) \cong \Z/7\Z \times \Z/(7 \cdot 23)\Z,
\]
and
\[
J(\F_{17}) \cong \Z/2\Z \times \Z/(2^2 \cdot 7^3 \cdot 31 \cdot 271)\Z .
\]
Hence $J(\Q)_{\mathrm{tors}}$
is isomorphic to a subgroup of $\Z/7\Z$.
Recall that the class $[c_0-c_\infty]$ has order $7$.
%As $X$ is not trigonal, and $c_0$, $c_\infty$
%have degree $3$, we see that $c_0$ and $c_\infty$ are not linearly
%equivalent, and so $[c_0-c_\infty]$ has exact order $7$.
Thus $[c_0-c_\infty]$ generates
$J(\Q)_{\mathrm{tors}}$.
Now
since $w_5$ interchanges $c_0$ and $c_\infty$,
\[
(w_5-1)([3c_0-3c_\infty])=6[c_\infty-c_0]=[c_0-c_\infty].
\]
Therefore $[c_0-c_\infty] \in A(\Q)$. We have now proved the following.
\begin{lem}\label{lem:A}
$A(\Q)=(\Z/7\Z) \cdot [c_0-c_\infty]$.
\end{lem}

\section{Proof of Theorem~\ref{thm:cubic_pts}}
In this section we prove Theorem~\ref{thm:cubic_pts} thereby completing
the proof of Theorem~\ref{thm:main}. Recall $X=X(\mathrm{b}5,\mathrm{ns}7)$.
Write $X^{(3)}$ for the $3$-rd symmetric power of $X$.
We shall prove the following result which immediately implies Theorem~\ref{thm:cubic_pts}.
\begin{prop}
$X^{(3)}(\Q)=\{c_0,\, c_\infty\}$.
\end{prop}
\begin{proof}
Let $x \in X^{(3)}(\Q)$.
%Then $x_{\F_3}=c_{0,{\F_3}}$ or $c_{\infty,\F_3}$.
%Let $x \in X^{(3)}(\Q)$.
By Lemma~\ref{lem:A} we have  $(1-w_5)[x-c_\infty]=\ell \cdot [c_0 - c_\infty]$
for some $\ell \in \Z/7\Z$. As $w_5(c_\infty)=c_0$ we may rewrite this as
\[
(x-w_5 (x)) \sim k \cdot (c_0-c_\infty)
\]
for some $k \in \{-3,\dotsc,3\}$.
We write $x_{\F_{3}}$, $c_{0,\F_3}$, $c_{\infty,\F_3} \in X^{(3)}(\F_3)$ for the reductions of $x$, $c_0$, $c_\infty$ modulo $3$
respectively.
It follows that
\begin{equation}\label{eqn:mwsieve}
\left(y- w_5(y)\right) \sim k \cdot \left( c_{0,{\F_3}} - c_{\infty,\F_3} \right)
\end{equation}
where $y=x_{\F_{3}}$. Using our model $D$ we enumerated $X^{(3)}(\F_3)$;
this has precisely $40$ elements.
For each $y \in X^{(3)}(\F_3)$ and for each $k=-3,\dotsc,3$
 we tested the relation \eqref{eqn:mwsieve} and found that it holds
only for $y=c_{0,\F_3}$ and $k=1$ and for $y=c_{\infty,\F_3}$
and $k=-1$. We therefore deduce that $x_{\F_3}=c_{0,\F_3}$ or $c_{\infty,\F_3}$.
We would like to conclude that $x=c_0$ or $c_\infty$. As $w_5$ swaps $c_0$ and $c_\infty$ and also their mod $3$ reductions,
we may suppose that $x_{\F_3}=c_\infty$. Let $\mu : X^{(3)} \rightarrow J$ be given by $z \mapsto [z-c_\infty]$
and $t: J \rightarrow A$ be simply $t=w_5-1$. Since $x_{\F_3}=c_{\infty,\F_3}$, the point $(t \circ \mu)(x) \in A(\Q)$
belongs to the kernel of reduction $A(\Q) \rightarrow A(\F_3)$. However as $A(\Q)$ is torsion,
this kernel of reduction is trivial \cite[Appendix]{Katz}. Thus $(t \circ \mu)(x)=0$.
To conclude that $x=c_\infty$ it is now enough to check that $t \circ \mu$ is a formal immersion at $c_{\infty,\F_3}$,
and for this we shall use the formal immersion criterion due to Derickx, Kamienny, Stein and Stoll \cite[Proposition 3.7]{DKSS}.

Write $\Omega_X \cong \Omega_J$ for the $6$-dimensional space of $1$-forms on $X/\F_3$.
We would like to write down the $4$-dimensional subspace $t^* (\Omega_A)$.
We easily do this since it is precisely that $-1$-eigenspace of $w_5^*$ acting on $\Omega_X$,
and we know the action of $w_5$ on our model $D$ from which can write down the corresponding
action on the $1$-forms. Let $\omega_1,\dotsc,\omega_4$
be an $\F_3$-basis for $t^*(\Omega_A)$. To check the formal immersion criterion of Derickx et al.\
at $c_{\infty,\F_3}$ we need to check that a certain $4 \times 3$ matrix
defined
in \cite[Proposition 3.7]{DKSS},
which we denote by $M$, has rank $3$.
As $3$ is inert in $\Q(\zeta_7)^+$,
we have
$c_{\infty,\F_3}=P_1+P_2+P_3$, where $P_i \in X(\F_{27})$ are distinct. This slightly
simplifies the description of the matrix $M$. Let $u_j \in \F_{27}(X)$
be a uniformizing element for $P_j$. Then $\omega_i/d{u_j}$ is a regular function
at $P_j$ and we may evaluate $(\omega_i/d{u_j}) (P_j) \in \F_{27}$. The matrix
is simply
\[
M=\left((\omega_i/d{u_j}) (P_j) \right)_{i=1,2,3,4; \, j=1,2,3}.
\]
We computed $M$ and checked that it has rank $3$ as required. This completes the proof.
\end{proof}

\end{document}